\documentclass[11pt]{article}
\usepackage{graphicx} 
\usepackage{url}
\usepackage{enumitem}
\usepackage{amsmath}
\usepackage{amssymb}
\usepackage{amsfonts}
\usepackage{breakcites}
\usepackage[fleqn]{nccmath}
\usepackage{ragged2e}
\usepackage{cite}

\usepackage{mathtools}
\newtheorem{Thm}{Theorem}
\newtheorem{Lem}{Lemma}
\newtheorem{Prop}{Proposition}

\renewcommand{\bf}[1]{\textbf{#1}}
\renewcommand{\it}[1]{\textit{#1}}
\newcommand{\equa}{\begin{equation}}
\newcommand{\equb}{\end{equation}}
\newcommand{\vs}{\vspace{+5mm}}
\newcommand{\eqa}{\begin{equation*}}
\newcommand{\eqb}{\begin{equation*}}
\newcommand{\eqn}{\begin{eqnarray*}}
\newcommand{\eqnb}{\end{eqnarray*}}

\newcommand{\bb}[1]{\mathbb{#1}}
\usepackage{lmodern} 
\usepackage{graphicx,scalerel}

\begin{document}

\begin{center}
    \LARGE{\bf{Long Time Behavior of General Markov Additive Processes}}\\
    \vs
    \normalsize{Celal Umut Yaran and Mine Çağlar}\\ 
    \normalsize{Koç University}
\end{center}
\begin{abstract}
\noindent
    We study general Markov additive processes when the state space of the modulator is a Polish space. Under some regularity assumptions, our main result is the characterization of the long-time behavior of the ordinate in terms of the associated ladder time process and the excursion measure. An important application of Markov additive processes is the Lamperti-Kiu transform, which gives a correspondence between $\bb{R}^d\backslash \{0\}$-valued self-similar Markov processes and $S^{d-1}\times \bb{R}$-valued Markov additive processes. The asymptotic behavior of the radial distance from the origin of a self-similar Markov process can be characterized by the long-time behavior of the ordinate of the corresponding Markov additive process. We show the applicability of our assumptions on some well-known self-similar Markov processes.\\
    \bf{Keywords:} Markov additive process, Fluctuation theory, Lamperti-Kiu transform, self-similar Markov process
\end{abstract}

\section{Introduction}
\flushleft
\justifying
In 1970s, Markov additive processes (MAP) were introduced in \cite{ccinlar1971markov,ccinlar1972markov}, essentially as a real-valued stochastic process $(\xi_t)_{t\geq 0}$ such that the increments of $\xi$ are modulated by a Markov process $(\Theta_t)_{t\geq 0}$. Formally, the process $(\xi,\Theta)$ is called a MAP if given $\{(\xi_s,\Theta_s),s\leq t\}$ for any $t\geq 0$, the process $(\xi_{s+t}-\xi_t,\Theta_{s+t})_{s\geq 0}$ has the same law as $(\xi_s,\Theta_s)_{s\geq 0}$ under $\bb{P}_{0,v}$ with $v=\Theta_t$, where $\bb{P}_{x,\theta}(\xi_0=x,\Theta_0=\theta)=1$. The following studies on general MAPs can be found first in \cite{ccinlar1973levy,ccinlar1976entrance,kaspi1982symmetric,kaspi1983excursions}, and in \cite{ney1987markov,neym1987markov}. Although MAPs have been introduced in such generality, most of the previous work has focused on the case the \it{modulator} $\Theta$ is a Markov chain with a countable state space (e.g. \cite{alsmeyer2018fluctuation,asmussen2010ruin,ccinlar1977shock}).

A key application of MAPs is the representation of real-valued self-similar Markov processes, which can be seen as a generalization of the relation between positive self-similar Markov processes and potentially killed Lévy processes given in \cite{lamperti1972semi}. A real-valued Feller process $(Z_t)_{t\geq 0}$ is called a real self-similar Markov process with index $\alpha>0$ if the process $(Z_t)_{t\geq 0}$ is equal in distribution to $(cZ_{c^{-\alpha}t})_{t\geq 0}$ for every $c>0$, with $P_z(Z_0=z)=1$. There exists a one-to-one correspondence between real self-similar Markov processes and $\bb{R}\times \{\pm 1 \}$-valued MAPs \cite{chaumont2013lamperti,kuznetsov2014hitting}. Here, the modulator $\Theta$ controls the sign of $Z$, and $e^\xi$ specifies how far away $Z$ is from the origin at time $t\geq 0$. A natural question whether a similar relation can be found for $\bb{R}^d \backslash \{ 0 \}$-valued self-similar Markov processes has been resolved by recognizing the radial distance from the origin of the process as a real-valued stochastic process modulated by the angular part of the process. There exists a one-to-one correspondence between $\bb{R}\times S^{d-1}$-valued MAPs and $\bb{R}^d\backslash \{ 0 \}$-valued self-similar Markov processes killed at the origin \cite{kiu1980semi}. This relation, known as Lamperti-Kiu transform, has brought back attention to MAPs when the state space of the modulator is a Polish space \cite{alili2017inversion,kyprianou2018stable,kyprianou2022stable,kyprianou2020entrance,kyprianou2023williams}. 

Let $S$ be a Polish space and $(\xi,\Theta)$ be an $\bb{R}\times S$-valued MAP with an invariant probability measure $\pi$ for the modulator $\Theta$. We are interested in the long-time behavior of the \it{ordinate} $\xi$. It has recently been shown in \cite{kyprianou2020entrance} that  $\xi_t$ exhibits exactly one of the three behaviors, namely, drifts to $+\infty$, oscillates, or drifts to $-\infty$ as $t\rightarrow \infty$, if the MAP satisfies a Harris-type condition and $\Theta$ has an invariant distribution. The proof follows from the strong law of large numbers for MAPs \cite[Prop.2.15]{kyprianou2020entrance} as inspired by a similar approach for Lévy processes. In fact, a MAP can be considered as a natural extension of a Lévy process. Therefore, further ideas from the fluctuation theory of Lévy processes can also be applicable to general MAPs\cite{bertoin1996levy}. In this paper, we use the distributions of last time at supremum $\bar{g}_\infty$  and infimum $\underline{g}_\infty$ arising from the corresponding ladder processes in the characterization of the long-time behavior of $\xi$. \\
\indent In Proposition \ref{Main Proposition}, we show under some regularity assumptions for  $(\xi,\Theta)$ that $\xi_t$ drifts to $-\infty$, that is, $\lim_{t\to\infty}\xi_t = -\infty$ $\:\bb{P}_{0,\pi}$-a.s., if \begin{equation*}
    \lim_{\lambda\to 0^+}\bb{E}_{0,\pi}[\exp(-\lambda \bar{g}_\infty)]=1 \quad \text{and} \quad \lim_{\lambda\to 0^+}\bb{E}_{0,\pi}[\exp(-\lambda \underline{g}_\infty)]=0,
\end{equation*}
as well as the other cases which imply that either the process $\xi_t$ oscillates, that is, \linebreak $\limsup_{t\to\infty}\xi_t = -\liminf_{t\to\infty}\xi_t = \infty$ or $\xi_t$ drifts to $+\infty$, that is, $\lim_{t\to\infty}\xi_t=\infty$. Indeed, for a real-valued process, $\lim_{\lambda \to 0^+}\bb{E}[\exp(-\lambda \bar{g}_\infty)]=1$ implies that $\bb{P}(\bar{g}_\infty < \infty )=1$ and $\lim_{\lambda \to 0^+}\bb{E}[\exp(-\lambda \underline{g}_\infty)]=0$ implies that $\bb{P}(\underline{g}_\infty <\infty )=0$. When these limits hold, the asymptotic behavior of a Lévy process $X$ is found as $\lim_{t\to\infty}X_t=-\infty$ in \cite{bertoin1996levy}. Along these lines, we determine the long-time behavior of $\xi$ depending on whether $\Bar{g}_\infty$ and $\underline{g}_\infty$ are finite or infinite using the strong Markov property and the invariant probability measure.

The main result of this paper is the characterization of the long-time behavior of the ordinate $\xi$ in terms of the associated ladder time process and the excursion measure as given in Theorem \ref{Main Theorem}.
The distribution of the last time at supremum at an exponentially distributed random time $e_q$ with parameter $q>0$, which is independent of the MAP, is given in \cite[Prop.2.3]{kyprianou2020entrance} as
\begin{equation*}
    \bb{E}_{0,\pi}[\exp(-\lambda \bar{g}_{e_q})] =\int_{\bb{R}^+ \times S \times \bb{R}^+}V_\pi^+ (dr,dv,dz) \exp(-\lambda r -qr) [q\ell^+(v) + \eta_v^+(\bf{1}_{e_q< \zeta^+})]
\end{equation*}
where $V^+_\pi$ is the potential measure associated with the corresponding ascending ladder process $(\Bar{L}_t^{-1},\xi^+_t,\Theta^+_t)$, $\eta_v^+$ denotes the excursion measure at the supremum, $\zeta^+$ denotes the lifetime of an excursion, and $\ell^+$ is a function related to the local time at the supremum $\Bar{L_t}$. Under a certain duality assumption called the reversibility property, which implies that $\underline{g}_t$ and $t-\bar{g}_t$ are equal in distribution\cite{alili2017inversion}, we first obtain the distribution of $\underline{g}_{e_q}$ in terms of the ascending ladder process. Then, based on this result and Proposition \ref{Main Proposition}, we show that $\xi_t$ drifts to $-\infty$ $\bb{P}_{0,\pi}$-a.s. if\begin{equation*}
    \int_{S\times \bb{R}^+}U^+_\pi (dv,dz)\eta^+_v(\zeta=\infty)=1
\end{equation*}
where \begin{equation*}
U^+_\pi (dv,dz)= \bb{E}_{0,\pi}\left[ \int_0^{\Bar{L}_\infty} \bf{1}_{\{ \Theta^+_s\in dv, \xi^+_s\in dz\}}ds \right].
\end{equation*} The conditions which determine either $\xi$ drifts to $+\infty$ or it oscillates are also examined in Theorem \ref{Main Theorem}. Furthermore, we extend this result by replacing the reversibility property with its weaker version, called the weak reversibility property \cite{kyprianou2020entrance}. In Theorem \ref{Generalised Theorem}, we examine the long-time behavior of the ordinate under the assumption of weak reversibility.

We study two important examples of well-known self-similar Markov processes to show that the corresponding MAPs via Lamperti-Kiu transform satisfy the assumptions of Theorem \ref{Main Theorem}. The first example is the free 2-dimensional Bessel process $X=(X^{(1)},X^{(2)})$ whose components $X^{(i)}$ are independent Bessel processes. The second example is the radial Dunkl process $X^W$ in $\bb{R}^2\backslash\{ 0 \}$. The long-time behavior of the ordinate of the MAP determines the asymptotic behavior of the radial distance from the origin of the corresponding $\bb{R}^d\backslash \{ 0 \}$-valued self-similar Markov process. Although the latter asymptotic behavior for both examples is known, we still 
prove some new properties of these processes using the Lamperti-Kiu transform. For instance, the modulators are Harris recurrent for both examples. We take $d=2$ for simplicity. However, we believe that these results can be generalized to higher dimensions. Moreover, the assumptions of Theorem \ref{Main Theorem} could be used for other classes of self-similar Markov processes with unknown asymptotic behavior. Further possible examples can be found in \cite{alili2019space}.

The paper is organized as follows. The fluctuation identities and duality properties of general MAPs are reviewed in Section 2. Our results on the long-time behavior of the ordinate of MAPs are presented in Section 3. In Section 4, we discuss some examples of self-similar Markov processes and their Lamperti-Kiu transforms in detail.

\section{Preliminaries}
\justifying
Let $(\xi_t,\Theta_t)_{t\geq 0}$ be a possibly killed Markov process on a filtered probability space \linebreak 
$(\Omega,\mathcal{F},(\mathcal{F}_t)_{t\geq 0},\bb{P}_{x,\theta})$ with $\mathbb{P}_{x,\theta}(\xi_0=x,\Theta_0=\theta)=1$ and state space $(\bb{R}\times S,\mathcal{B}(\bb{R}\times S))$ with an extra isolated state $\partial$, where $S$ is a Polish space, $\mathcal{B}$ denotes the Borel $\sigma$-algebra, and $(\mathcal{F}_t)_{t\geq 0}$ is the minimal augmented admissable filtration.
The process $(\xi_t,\Theta_t)_{t\geq 0}$ is called a MAP on $\bb{R}\times S$, if \begin{equation*}
     \bb{P}_{0,\theta}( (\xi_{t+s}-\xi_t)\in \Gamma \, , \, \Theta_{t+s}\in A \, | \, \mathcal{F}_t ) = \bb{P}_{0,v}(\xi_s \in \Gamma \, , \, \Theta_s\in A)
    \end{equation*}
for all $\theta\in S$, $\Gamma\in \mathcal{B}( \bb{R} )$, $A \in \mathcal{B} (S)$, and $t\geq 0$,  where $v=\Theta_t$, and denoted by $((\xi,\Theta),\bb{P})=((\xi_t,\Theta_t)_{t\geq 0},\bb{P}_{0,\theta})$. Let $\zeta:=\inf \{ t>0 : (\xi_t,\Theta_t)=\partial \}$ denote the lifetime of the MAP. We call $\xi$ as the ordinate and $\Theta$ as the modulator of the MAP. By the defining property, a MAP is translation invariant in $\xi$, that is, $(\xi_t,\Theta_t)$ under $\bb{P}_{x,\theta}$ has the same law as $(\xi_t+x,\Theta_t)$ under $\bb{P}_{0,\theta}$ for all $x\in \bb{R}$ and $\theta\in S$. A MAP has strong Markov property, that is, for a stopping time $T$, we have
\begin{equation}\label{StrongMAP}
     \bb{P}_{0,\theta}( (\xi_{T+t}-\xi_T)\in \Gamma \, , \, \Theta_{T+t}\in A \, | \, \mathcal{F}_T ) = \bb{P}_{0,v}(\xi_t \in \Gamma \, , \, \Theta_t\in A)
    \end{equation}
for all $\theta\in S$, $\Gamma\in \mathcal{B}( \bb{R} )$, $A \in \mathcal{B} (S)$, and $t>0$, where $v=\Theta_T$ \cite[Thm.3.2]{ccinlar1972markov}.

    Throughout this section, we assume that $\Theta$ is a Hunt process, $\xi$ is quasi-left continuous on $[0,\zeta]$ where $\zeta$ is the lifetime of the process and the MAP $(\xi,\Theta)$ is upwards regular, that is, \begin{equation*}
        \bb{P}_{0,\theta}(\tau_0^+ = 0 ) = 1 
    \end{equation*}
    $\forall \theta \in S$ where $\tau^+_0 := \inf \{ t>0 \, : \, \xi_t>0\}$. Let $\Bar{\xi}_t:=\text{sup}_{s\leq t}\, \xi_s$ denote the running supremum process for the ordinate $\xi$ and  $U_t:=\Bar{\xi}_t-\xi_t$ be the corresponding reflected process. Denote the zeros of the reflected process by $\Bar{M}:=\{ t\geq 0 \, : \, U_t=0\}$ and let $\Bar{M}^{cl}$ be the closure of $\Bar{M}$ in $\bb{R}^+$. The complement of $\Bar{M}^{cl}$ consists of open random intervals. For an interval with endpoints $g$ and $d$, define\begin{equation*}
        (\epsilon_s^{(g)},v_s^{(g)}):=
 \begin{cases}
        (U_{g+s},\Theta_{g+s}) & \text{if } 0\leq g+s\leq d\\
        (U_d,\Theta_d) &  \text{if } d\leq g+s.
    \end{cases}
    \end{equation*}
    Notice that $(\epsilon_s^{(g)},v_s^{(g)})$ corresponds to an excursion from the supremum. Let $\zeta^+$ denote the lifetime of the excursion and $\eta_\theta^+$ be the excursion measure defined in \cite[\S 2.2]{kyprianou2020entrance}.
    Let $\Bar{L}$ be the local time at the supremum associated to the ordinate $\xi$. Then, by \cite{kaspi1983excursions}, there exists a measurable function $\ell^+:S\to \bb{R}^+$ such that\begin{ceqn}\begin{equation*}
        \int_0^t \bf{1}_{\{ s\in \Bar{M}\}}ds=\int_0^t \bf{1}_{\{ s\in \Bar{M}^{cl}\}}ds=\int_0^t \ell^+(\Theta_s)d\Bar{L}_s \quad\quad \bb{P}_{0,\pi}\text{-a.s.}
    \end{equation*}\end{ceqn}
    Let $\Bar{L}_t^{-1}:= \text{inf}\{ s\, : \, L_s>t\}$ be the inverse local time process and define $\xi_t^+=\xi_{\Bar{L}_t^{-1}}$ and $\Theta_t^+=\Theta_{\Bar{L}_t^{-1}}$ for all $t$ such that $\Bar{L}_t^{-1}<\infty$, and otherwise assign both $\xi_t^+$ and $\Theta_t^+$ to the cemetery state $\partial$. Then, $(\Bar{L}_t^{-1},\xi_t^+,\Theta_t^+)_{t\geq 0}$ defines a MAP called \it{ascending ladder process}, where both $\Bar{L}_t$ and $\xi^+_t$ are ordinates \cite{kyprianou2020entrance}. Let $\Bar{g}_t=\text{sup}\{ s\leq t \, : \, s\in \Bar{M}^{cl}\}$ and $\Bar{\Theta}_t= \Theta_{\Bar{g}_t} \bf{1}_{\{ \Bar{\xi}_t=\xi_{\Bar{g}_t}\} } + \Theta_{\Bar{g}_t^-}\bf{1}_{\{ \Bar{\xi}_t >\xi_{\Bar{g}_t} \} }$ and recall that $e_q$ is an exponentially distributed random time with parameter $q>0$, which is independent of $(\xi,\Theta)$. Then, by \cite[Prop.2.3]{kyprianou2020entrance}, for bounded measurable functions $F,G:[0,\infty)\times \bb{R}\times S \to \bb{R}$ and for $\theta\in S$, we have
\begin{multline}\label{Kyprianou Proposition}
\bb{E}_{0,\theta}\left[G(\Bar{g}_{e_q},\Bar{\xi}_{e_q},\Bar{\Theta}_{e_q})F(e_q-\Bar{g}_{e_q},\Bar{\xi}_{e_q}-\xi_{e_q},\Theta_{e_q})   \right]=\\ 
\int_{\bb{R}^+\times\bb{R}^+\times S}e^{-qr}G(r,z,v)\left[ q\ell^+(v)F(0,0,v)
+\eta^+_v(F(\tau,\epsilon_\tau,v_\tau)\bf{1}_{\{e_q<\zeta^+\}} \right] V_\theta^+(dr,dv,dz)
\end{multline}
where \begin{equation*}
    V_\theta^+(dr,dv,dz):=\bb{E}_{0,\theta}\left[ \int_0^{\Bar{L}_\infty} \bf{1}_{\{ \Bar{L}^{-1}_s\in dr , \Theta^+_s\in dv, \xi^+_s\in dz\}}ds \right].
\end{equation*}
\bf{Remark 1} We say that the MAP $((\xi,\Theta),\bb{P})$ is downwards regular if\begin{equation*}
    \bb{P}_{0,\theta}(\tau^-_0 = 0 ) =1
\end{equation*} 
$\forall \theta \in S$ where $\tau^-_0 := \inf \{ t>0 \, : \, \xi_t <0 \}$. Now, suppose that the MAP is downwards regular. Then, we can construct the \it{descending ladder process} $(\underline{L}_t^{-1},\xi^-_t,\Theta^-_t)$ by considering the ascending ladder process corresponding to the MAP $( (\xi_t,\Theta_t)_{t\geq 0}\, , \, \hat{\bb{P}}_{-x,\theta} )$ obtained from $( (\xi_t,\Theta_t)_{t\geq 0}\, , \, \bb{P}_{x,\theta} )$ by replacing $\xi$ by its negative $-\xi$. In that case, the running infimum $\underline{\xi}_t := \inf_{s\leq t}\xi_s$ under $\bb{P}_{0,\theta}$ corresponds to the running supremum $\overline{\xi}_t := \sup_{s\leq t}(\xi_s)$ under $\Hat{\bb{P}}_{0,\theta}$. Similarly, the last time at the infimum $\underline{g}_t$, the lifetime of the excursion from the infimum $\zeta^-$, the excursion measure at the the infimum $\eta^-$ and the local time at the infimum $\underline{L}_t$ with density $\ell^-$ under $\bb{P}_{x,\theta}$ correspond to ${\Bar{g}}_t$, $\zeta^+$, ${\eta}^+$, ${\Bar{L}}_t$ and ${\ell}^+$ under $\Hat{\bb{P}}_{-x,\theta}$, respectively. Using these notations, the excursion from the infimum version of Equation \eqref{Kyprianou Proposition} can be stated with potential measure associated with $(\underline{L}_t^{-1},\xi^-_t,\Theta^-_t) $ defined as\begin{eqnarray*}
    V_\theta^-(dr,dv,dz)&:=&\bb{E}_{0,\theta}\left[ \int_0^{\underline{L}_\infty} \bf{1}_{\{ \underline{L}^{-1}_s\in dr , \Theta^-_s\in dv, \xi^-_s\in dz\}}ds \right] \\
    & = & \Hat{\bb{E}}_{0,\theta}\left[ \int_0^{{\bar{L}}_\infty} \bf{1}_{\{ {\Bar{L}}^{-1}_s\in dr , \Theta^+_s\in dv, \xi_s^+\in dz\}}ds \right]
\end{eqnarray*} 
\cite[\S 2.5]{kyprianou2023williams}.

Let $E$ be a Polish space. A pair of $E$-valued Markov processes with respective transition probabilities $(P_t)_{t\geq 0}$ and $(\widetilde{P}_t)_{t\geq 0}$ are set to be \it{weak duality} with respect to some $\sigma$-finite measure $m(dx)$ if for all positive measurable functions $f$ and $g$ \begin{equation*}
    \int_E g(x)P_tf(x) m(dx) = \int_E f(x)\widetilde{P}_tg(x)m(dx).
\end{equation*}
The details of weak duality of Markov processes can be found in \cite[Ch.13]{chung2006markov}. Consider the MAPs $((\xi,\Theta),\bb{P})$ and $((\xi,\Theta),\Hat{\bb{P}})$. The problem is to determine whether these two MAPs are in weak duality with respect to the measure $m(ds)dy$ where $dy$ is the Lebesgue measure and $m(ds)$ is some $\sigma$-finite measure on $\mathcal{B}(S)$. For weak duality, \cite[Lem.3.1]{alili2017inversion} provides a practical criterion. Namely, the MAPs $((\xi,\Theta),\bb{P})$ and $( (\xi,\Theta) , \hat{\bb{P}})$ are in weak duality with respect to measure $m(ds)dy$ if and only if \begin{equation}\label{Reversibility}
    \bb{P}_{0,v}(\xi_t\in dx, \Theta_t\in dw)m(dv)=\bb{P}_{0,w}(\xi_t\in dx, \Theta_t\in dv)m(dw),
\end{equation}
which is called the \it{reversibility property} of $( (\xi,\Theta),\bb{P})$. We have a more general version of this property from \cite{kyprianou2020entrance}. Let $( (\xi,\Theta) ,  \tilde{\bb{P}} )$ be another MAP and this time denote by  $((\xi,\Theta) ,  \hat{\bb{P}} )$ the process obtained from $( (\xi,\Theta) , \tilde{\bb{P}} )$ by replacing $\xi$ by its negative. Assume that the MAPs $( (\xi,\Theta),\bb{P})$ and $( (\xi,\Theta),\tilde{\bb{P}})$ satisfy the relation \begin{equation}\label{WeakReversibility}
     \bb{P}_{0,v}(\xi_t\in dx, \Theta_t\in dw)m(dv)=\tilde{\bb{P}}_{0,w}(\xi_t\in dx, \Theta_t\in dv)m(dw)
\end{equation}
where $m(ds)$ is some $\sigma$-finite measure on $S$. Then, Equation (\ref{WeakReversibility}) is called the \it{weak reversibility property} of MAPs $((\xi,\Theta),\bb{P})$ and $( (\xi,\Theta) , \tilde{\bb{P}})$. Notice that the reversibility property is a particular case of Equation (\ref{WeakReversibility}) if $\tilde{\bb{P}}
= \bb{P}$. Under the weak reversibility property, we have the following duality properties borrowed from \cite[Lem.3.8,Prop.3.9]{kyprianou2020entrance}.
\begin{Lem}\label{Duality Lemma}
Consider a MAP $( (\xi,\Theta),\bb{P})$ such that there exists a MAP $( (\xi,\Theta) ,  \tilde{\bb{P}})$ with weak reversibility property. Let $\hat{\bb{P}}_{0,\theta}$ be the law of $(-\xi,\Theta)$ under $\tilde{\bb{P}}_{0,\theta}$. Then, for every $t>0$, we have \begin{enumerate}
    \item for $x\in\bb{R}$, the process 
    $(\xi_{(t-s)-}-\xi_t,\Theta_{(t-s)-})_{0\leq s \leq t}$  under $\bb{P}_{x,m}$ has the same law as $(\xi_s,\Theta_s)_{0\leq s \leq t}$ under $\hat{\bb{P}}_{0,m}$.
\item $(\Theta_0, t-\Bar{g}_t, \Theta_t, \Bar{\xi}_t-\xi_t, \Bar{g}_t, \Bar{\Theta}_t, \Bar{\xi}_t)$ under $\hat{\bb{P}}_{0,m}$ is equal in distribution to\\ $(\Theta_t,\Bar{g}_t,\Theta_0,\Bar{\xi}_t,t-\Bar{g}_t,\Bar{\Theta}_t,\Bar{\xi}_t-\xi_t)$ under $\bb{P}_{0,m}.$
\end{enumerate}
\end{Lem}
\section{Long Time Behavior}
It has recently been shown in \cite{kyprianou2020entrance} that  $\xi_t$ exhibits exactly one of the three behaviors, namely, drifts to $+\infty$, oscillates,  or drifts to $-\infty$ as $t\rightarrow \infty$ almost surely if the MAP satisfies a Harris-type condition and $\Theta$  has an invariant distribution. A regenerative structure defined by the Harris-type condition \cite{ney1987markov} is closely related to the ergodic properties of the process. Therefore, the strong law of large numbers is proved first as inspired by a similar approach for Lévy processes, and is given by $\lim_{t\rightarrow \infty} \xi_t/t = \mathbb{E}_{0,\pi}[\xi _1]$, $\mathbb{P}_{0, \theta}$-a.s. for every $\theta \in {\cal S}$. The trichotomy follows whether $\mathbb{E}[\xi _1]$ is positive, negative, or equal to zero \cite{kyprianou2020entrance}.

Another idea for determining the long-time behavior of the ordinate is the use of fluctuation theory, also borrowed from Lévy processes \cite{bertoin1996levy}. The distributions of the last time at supremum $\bar{g}_\infty$  and the infimum $\underline{g}_\infty$ arising from the corresponding ladder processes can characterize the long-time behavior of $\xi$. We follow this direction without assuming Harris-type recurrence, but only the existence of an invariant distribution $\pi$. Therefore, the characterization is a $\mathbb{P}_{0, \pi}$-a.s.result rather than holding for any starting $\theta \in {\cal S}$. We make the following assumptions in this section.\\
\bf{Assumption 1}: $\Theta$ is a Hunt process with invariant probability measure $\pi$ and $\xi$ is quasi-left continuous.\\
\bf{Assumption 2}: The MAP $((\xi,\Theta),\bb{P})$ is both upwards and downwards regular.\\
\bf{Assumption 3}: $\bb{E}_{0,\pi}[\sup_{s\in [0,1]}|\xi_s|]$ is finite.

The following proposition gives a criterion for the long-time behavior of MAPs in terms of the distributions of $\Bar{g}_t$ and $\underline{g}_t$. In \cite[Prop.2.16]{kyprianou2020entrance}, this behavior is determined in terms of $\bb{E}_{0,\pi}[\xi_t]$ with some extra assumptions using the strong law of large numbers. Therefore, we achieve a different type of characterization thanks to the alternative approach resulting from the generalization of \cite[Ch.6.3]{bertoin1996levy}.
\begin{Prop}\label{Main Proposition}
    Consider a MAP $((\xi,\Theta),\bb{P})$, such that Assumptions 1,2, and 3 are satisfied. Then, we have\begin{enumerate}[label=(\roman*)]
 \item $\xi_t$ drifts to $-\infty$ $\:\bb{P}_{0,\pi}$-a.s., if and only if
 $\Bar{g}_\infty < \infty$  and $\underline{g}_\infty = \infty$   $\:\bb{P}_{0,\pi}$-a.s.
 \item $\xi_t$ drifts to $+\infty$ $\bb{P}_{0,\pi}$-a.s., if and only if $\Bar{g}_\infty = \infty$  and $\underline{g}_\infty < \infty$  $\:\bb{P}_{0,\pi}$-a.s.
\item $\xi_t$ oscillates $\bb{P}_{0,\pi}$-a.s., if and only if $\Bar{g}_\infty = \infty$  and $\underline{g}_\infty = \infty$  $\:\bb{P}_{0,\pi}$-a.s.
\end{enumerate}
\end{Prop}
\bf{Proof}: The necessity part is clear, thus we proceed by showing the sufficiency part.\\
\bf{(i)} Since the MAP is downwards regular by Assumption 2, and $\underline{g}_\infty =\infty$ $\:\bb{P}_{0,\pi}$-a.s., it follows that, $\xi$ visits $(-\infty,-x]$ $\forall x>0$ and the passage time $T_{(-\infty,-x]}:=\inf \{ t>0 \, | \, \xi_t\in (-\infty,-x] \}$ is $\bb{P}_{0,\pi}$-a.s. finite.\\
Assumption 3 implies that $\bb{E}_{0,\pi} [\sup_{s\in [0,t]}| \xi_s | ]<\infty$ for all $t>0$ by \cite[Lem.2.13]{kyprianou2020entrance}. Hence,  $\Bar{\xi}_\infty<\infty$ $\:\bb{P}_{0,\pi}$-a.s. since $\Bar{g}_\infty<\infty$ $\:\bb{P}_{0,\pi}$-a.s. It means that, we can pick $x_\epsilon$ large enough such that\begin{equation*}
    \bb{P}_{0,\pi}(\xi_t>x/2\,\,\,\text{for some }\, t\geq 0)\leq \epsilon \quad \forall x\geq x_\epsilon.
\end{equation*}
Denote $T:=T_{(-\infty,-x]}$. Then, for $x\geq x_\epsilon$, we have
\begin{flalign}
         &   \bb{P}_{0,\pi}(\xi_t>-x/2\,\,\,\text{for some }\, t\geq T )& \nonumber\\
            =\,\, & \bb{P}_{0,\pi}(\xi_{t+T}>-x/2\,\,\,\text{for some }\, t\geq 0 )& \nonumber \\
                \leq \,\,& \bb{P}_{0,\pi}(\xi_{t+T(-\infty ,-x]}-\xi_{T(-\infty ,-x]}>x/2\,\,\,\text{for some }\, t\geq 0  ) & \label{less than -x} \\
                =\, & \int_{S} \bb{P}_{0,\theta}(\xi_{t+T}-\xi_T>x/2\,\text{for some }\, t\geq 0 ) \pi(d\theta) \label{first integral}&  \\
                =\, & \int_{S} \int_{S} \bb{P}_{0,\theta}(\xi_{t+T}-\xi_T>x/2\,\text{for some }\, t\geq 0\, | \, \Theta_T =\phi )\bb{P}_{0,\theta}(\Theta_T \in d\phi) \pi(d\theta)& \nonumber \\
                 =\, & \int_{S}\int_{S} \bb{P}_{0,\phi}(\xi_t>x/2\,\,\,\text{for some }\, t\geq 0)\,\bb{P}_{0,\theta}(\Theta_T\in d\phi )   \pi(d\theta )& \label{usage of strong Markov} \\
                =\, & \int_{S} \bb{P}_{0,\phi}(\xi_t>x/2\,\,\,\text{for some }\, t\geq 0) \,\bb{P}_{0,\pi}(\Theta_T \in d\phi) &\nonumber \\
                =\, & \int_{S} \bb{P}_{0,\phi}(\xi_t>x/2\,\,\,\text{for some }\, t\geq 0)\pi(d\phi)& \label{last integral} \\
                =\,\, & \bb{P}_{0,\pi}(\xi_t>x/2\,\,\,\text{for some }\, t\geq 0) \quad < \quad \epsilon \label{last step} &
                 \end{flalign}
where we use $\xi_{T(-\infty ,-x]} \leq -x$ $\:\bb{P}_{0,\pi}$-a.s. in (\ref{less than -x}), definition of $\bb{P}_{0,\pi}$ in (\ref{first integral}) and (\ref{last step}), strong Markov property (\ref{StrongMAP}) in (\ref{usage of strong Markov}) and the fact that $\pi$ is an invariant measure in (\ref{last integral}). Hence, $\xi$ drifts to $-\infty$ $\:\bb{P}_{0,\pi}$-a.s.\\
\bf{(ii)} Apply the previous part of the proof to the process $( (\xi,\Theta) ,  \Hat{\bb{P}})$ using the upwards regularity of $((\xi,\Theta),\bb{P})$, where $\Hat{\bb{P}}$ is defined by replacing $\xi$ by its negative. \\
\bf{(iii)} Since the MAP is both upwards and downwards regular by Assumption 2, $\xi$ visits $(-\infty,-x]$ and $[x,\infty)$ $\forall x>0$ and both $T_{(-\infty,-x]}$ and $T_{[x,\infty)}$ are finite $\bb{P}_{0,\pi}$-a.s. Hence, $\xi$ oscillates.  $\hfill \Box$\\
\bf{Assumption 4}: The reversibility property \eqref{Reversibility} holds.\\
\indent Consider a MAP $((\xi,\Theta),\bb{P})$ such that also Assumption 4 is satisfied. By taking $\tilde{\bb{P}}=\bb{P}$ in the second part of Lemma \ref{Duality Lemma}, thanks to the reversibility property \eqref{Reversibility}, we can deduce that $\underline{g}_t$ and $t-\bar{g}_t$ have the same distribution under $\bb{P}_{0,\pi}$ as they are both equal in distribution to ${\bar{g}}_t$ under $\hat{\bb{P}}_{0,\pi}$ where $\hat{\bb{P}}_{0,\pi}$ is the law of $(-\xi,\Theta)$ under $\bb{P}_{0,\pi}$. Hence, using Proposition \ref{Main Proposition}, we obtain the following identities
\begin{eqnarray}
    \bb{E}_{0,\pi}[e^{-\lambda \Bar{g}_{e_q}}]&=&\int_{\bb{R}^+\times S \times \bb{R}^+}V_\pi^+(dr,dv,dz)e^{-\lambda r -qr}[q\ell^+(v)+\eta_v^+(\bf{1}_{\{e_q<\zeta^+\}})] \label{last time at supremum} \\
    \bb{E}_{0,\pi}[e^{-\lambda \underline{g}_{e_q}}]&=&\int_{\bb{R}^+\times S \times \bb{R}^+}V_\pi^+(dr,dv,dz)e^{ -qr}[q\ell^+(v)+\eta_v^+(e^{-\lambda e_q} \bf{1}_{\{e_q<\zeta^+\}})]\label{last time at minimum}
\end{eqnarray}
where \begin{equation*}
    V_\pi^+(dr,dv,dz)= \bb{E}_{0,\pi}\left[ \int_0^{\Bar{L}_\infty} \bf{1}_{\{ \Bar{L}^{-1}_s\in dr \, , \, \Theta^+_s\in dv \, , \, \xi^+_s\in dz \} }ds \right].
\end{equation*}
 We have assumed that MAP is also downwards regular, so we can apply \eqref{Kyprianou Proposition} for the process $((\xi,\Theta), \hat{\bb{P}})$. We obtain the following analogs of (\ref{last time at supremum}) and (\ref{last time at minimum}) as\begin{eqnarray}
    \Hat{\bb{E}}_{0,\pi}[e^{-\lambda {\Bar{g}}_{e_q}}]&=&\int_{\bb{R}^+\times S \times \bb{R}^+}V_\pi^-(dr,dv,dz)e^{-\lambda r -qr}[q\ell^-(v)+\eta_v^-(\bf{1}_{\{e_q<\zeta^-\}})]\label{last time at supremum of dual}
\\
        \Hat{\bb{E}}_{0,\pi}[e^{-\lambda {\underline{g}}_{e_q}}]&=&\int_{\bb{R}^+\times S \times \bb{R}^+}V_\pi^-(dr,dv,dz)e^{ -qr}[q\ell^-(v)+\eta_v^-(e^{-\lambda e_q} \bf{1}_{\{e_q<\zeta^-\}})]\nonumber
\end{eqnarray}
where \begin{equation*}
    V_\pi^-(dr,dv,dz)= \bb{E}_{0,\pi}\left[ \int_0^{\underline{L}_\infty} \bf{1}_{\{ \underline{L}^{-1}_s\in dr \, , \, \Theta^-_s\in dv \, , \, \xi^-_s\in dz \} }ds \right].
\end{equation*}
\bf{Remark 2} If we only use identities (\ref{last time at supremum}) and (\ref{last time at minimum}), we will not be able to distinguish the cases where $\xi_t$ drifts to $+\infty$ or oscillates. Therefore, we have also obtained these identities in terms of the descending ladder process $(\underline{L}_t^{-1},\xi^-_t,\Theta^-_t)$.

The following is our main theorem.
\begin{Thm}\label{Main Theorem}
    Let $(\xi,\Theta),\bb{P})$ be a MAP such that Assumptions 1,2,3, and 4 are satisfied. Then, we have
\begin{enumerate}[label=(\roman*)]
    \item $\xi_t$ drifts to $-\infty$ $\:\bb{P}_{0,\pi}$-a.s., if and only if\begin{equation}\label{MainThm(i)}
    \int_{S\times \bb{R}^+}U_\pi^+(dv,dz)\eta^+_v(\zeta^+ =+\infty)=1
\end{equation}
\item  $\xi_t$ drifts to $+\infty$ $\:\bb{P}_{0,\pi}$-a.s.,  if and only if\begin{equation}\label{MainThm(ii)}
    \int_{S\times \bb{R}^-}U_\pi^-(dv,dz)\eta^-_v(\zeta^- =+\infty)=1
\end{equation}
\item $\xi_t$ oscillates $\bb{P}_{0,\pi}$-a.s., if and only if \begin{equation}\label{MainThm(iii)}
    \int_{S\times \bb{R}^+}U_\pi^+(dv,dz)\eta^+_v(\zeta^+ =+\infty)=0 \quad \text{and} \quad \int_{S\times \bb{R}^+}U_\pi^-(dv,dz)\eta^-_v(\zeta^- =+\infty)=0
    \end{equation}
\end{enumerate}
where \begin{equation*}
U^+_\theta (dv,dz)= \bb{E}_{0,\theta}\left[ \int_0^{\Bar{L}_\infty} \bf{1}_{\{ \Theta^+_s\in dv, \xi^+_s\in dz\}}ds \right] \quad \text{and} \quad
     U^-_\theta (dv,dz)= 
     \bb{E}_{0,\pi}\left[ \int_0^{\underline{L}_\infty} \bf{1}_{\{ \Theta^-_s\in dv \, , \, \xi^-_s\in dz \} }ds \right].
\end{equation*}
\end{Thm}
\bf{Proof}: \bf{(i)} We can rewrite (\ref{last time at supremum}) as \begin{equation*}
    \bb{E}_{0,\pi}\left[e^{-\lambda \Bar{g}_{e_q}}\right]=\int_{\bb{R}^+\times S \times \bb{R}^+}V_\pi^+(dr,dv,dz)e^{-\lambda r -qr}[q\ell^+(v)+\eta_v^+(1-e^{-q\zeta^+})].
\end{equation*} 
For any $0<q<\lambda/2$, the integrand is bounded from above by\begin{equation*}
    e^{-\lambda r}\left[\frac{\lambda}{2}\ell^+(v) + \eta_v^+(1-e^{-\lambda \zeta^+/2})\right].
\end{equation*}
Also, we have\begin{equation*}
    \int_{\bb{R}^+\times S \times \bb{R}^+}V_\pi^+(dr,dv,dz)e^{-\lambda r}\left[\frac{\lambda}{2}\ell^+(v)+\eta_v^+(1-e^{-\lambda \zeta^+/2})\right]=\bb{E}_{0,\pi}\left[e^{-\frac{\lambda}{2} \Bar{g}_{e_{\lambda/2}}}\right]<+\infty.
\end{equation*}
Then, by applying the dominated convergence theorem on both sides of (\ref{last time at supremum}) as $q\to 0^+$, we get\begin{equation*}
    \bb{E}_{0,\pi}\left[e^{-\lambda \Bar{g}_\infty}\right]=\int_{\bb{R}^+\times S \times \bb{R}^+}V_\pi^+(dr,dv,dz)e^{-\lambda r }\eta_v^+(\zeta^+=+\infty).
\end{equation*}
By letting $\lambda\to 0^+$, we obtain\begin{equation}\label{integral1}
    \lim_{\lambda \to 0^+} \bb{E}_{0,\pi}\left[e^{-\lambda \Bar{g}_\infty}\right]=\int_{S \times \bb{R}^+}U_\pi^+(dv,dz)\eta_v^+(\zeta^+=+\infty).
\end{equation}
Note that if $\xi_t$ drifts to $-\infty$ $\bb{P}_{0,\pi}$-a.s., then the left-hand side of \eqref{integral1} is equal to $1$ by the first part of Proposition \ref{Main Proposition}. Hence \eqref{MainThm(i)} holds.\\
Now, assume that \eqref{MainThm(i)} holds. Then, \eqref{integral1} implies that $\lim_{\lambda \to 0^+} \bb{E}_{0,\pi}\left[e^{-\lambda \Bar{g}_\infty}\right]=1$.
Hence, we have shown that $\Bar{g}_\infty<\infty$ $\:\bb{P}_{0,\pi}$-a.s. By considering (\ref{last time at minimum}), we have\begin{equation}\label{exact form of last time at minimum}
     \bb{E}_{0,\pi}\left[e^{-\lambda \underline{g}_{e_q}}\right]=\int_{\bb{R}^+\times S \times \bb{R}^+}V_\pi^+(dr,dv,dz)e^{ -qr}\left[ q\ell^+(v)+\frac{q}{q+\lambda}\eta_v^+(1-e^{-(\lambda+q)\zeta^+})\right].
\end{equation}
Since $\Bar{g}_\infty$ is $\bb{P}_{0,\pi}$-a.s. finite, for each $0<\epsilon <1$, there exists $t_\epsilon>0$ such that\begin{equation*}
    \bb{P}_{0,\pi}( \Bar{g}_\infty>t_\epsilon )< \epsilon
\end{equation*}
which implies that\begin{equation*}
     \bb{P}_{0,\pi}(\Bar{L}_\infty>t_\epsilon )< \epsilon .
\end{equation*}
We choose $(t_\epsilon)_{\epsilon \geq 0}$ such that $t_\epsilon$ increases to infinity as $\epsilon$ goes to zero. Note that, for each $\epsilon>0$, we can rewrite the right-hand side of (\ref{exact form of last time at minimum}) as\begin{multline}\label{two part of integral}
   \int_{[0,t_\epsilon]\times S \times \bb{R}^+}V_\pi^+(dr,dv,dz)e^{ -qr}\left[ q\ell^+(v) +\frac{q}{q+\lambda}\eta_v^+(1-e^{-(\lambda+q)\zeta^+})\right] \\ +
     \int_{(t_\epsilon,\infty)\times S \times \bb{R}^+}V_\pi^+(dr,dv,dz)e^{ -qr}\left[ q\ell^+(v)+\frac{q}{q+\lambda}\eta_v^+(1-e^{-(\lambda+q)\zeta^+})\right].
\end{multline}
The second integral in (\ref{two part of integral}) can be rewritten as\begin{equation*}
      \bb{E}_{0,\pi}\left[ \int_{\Bar{L}_{t_\epsilon}}^{\Bar{L}_\infty} ds\, e^{-q\Bar{L}_s^{-1}} \left[ q\ell^+(\Theta^+_s) + \frac{q}{q+\lambda} \eta^+_{\Theta_s^+}(1-e^{-(\lambda+q)\zeta^+})    \right]      \right].
\end{equation*}
Note that as $\epsilon$ goes to $0$,  $t_\epsilon$ increases to infinity. It follows that $\Bar{L}_{t_\epsilon}$ increases to $\Bar{L}_\infty$. Furthermore, the integrand is positive, so the integral decreases monotonically as $\epsilon$ goes to $0$. Hence, by the monotone convergence theorem, we have\begin{multline*}
    \lim_{\epsilon\to 0} \bb{E}_{0,\pi}\left[ \int_{\Bar{L}_{t_\epsilon}}^{\Bar{L}_\infty} ds\, e^{-q\Bar{L}_s^{-1}} \left[ q\ell^+(\Theta^+_s) + \frac{q}{q+\lambda} \eta^+_{\Theta_s^+}(1-e^{-(\lambda+q)\zeta^+})    \right]    \right]=\\
    \bb{E}_{0,\pi} \left[\lim_{\epsilon \to 0} \int_{\Bar{L}_{t_\epsilon}}^{\Bar{L}_\infty} ds\, e^{-q\Bar{L}_s^{-1}} \left[ q\ell^+(\Theta^+_s) + \frac{q}{q+\lambda} \eta^+_{\Theta_s^+}(1-e^{-(\lambda+q)\zeta^+})    \right]    \right]=0.
    \end{multline*}
It means that we can fix a sufficiently small $\epsilon>0$ such that the second integral in (\ref{two part of integral}) is bounded.
For fixed $\lambda$ and $\epsilon >0$, there exists an integer $n\in \bb{N}$ such that $\frac{-\lambda t_\epsilon +\sqrt{\lambda^2 t_\epsilon^2+4\lambda t_\epsilon}}{2t_\epsilon}>\frac{\lambda}{n}$. Then, for any $0<q<\lambda/n$, the integrand in the first integral in (\ref{two part of integral}) is bounded from above by\begin{equation*}
    e^{ -\frac{\lambda}{n}r}\left[\frac{\lambda}{n}\ell^+(v)+\frac{1}{n+1}\eta_v^+(1-e^{-\frac{n+1}{n}\lambda\zeta^+})\right].
\end{equation*}
Therefore, the first integral in \eqref{two part of integral} is bounded from above by\begin{equation*}
  \int_{\bb{R}^+\times S \times \bb{R}^+}V_\pi^+(dr,dv,dz)e^{ -\frac{\lambda}{n}r}\left[\frac{\lambda}{n}\ell^+(v)+\frac{1}{n+1}\eta_v^+(1-e^{-\frac{n+1}{n}\lambda\zeta^+})\right] 
    = \bb{E}_{0,\pi}\left[e^{-(\lambda/n) \underline{g}_{e_{\lambda/n}}} \right]<\infty.
\end{equation*}
Then, applying the dominated convergence theorem as $q\to 0^+$ on both sides of (\ref{last time at minimum}) yields
\begin{equation*}
    \bb{E}_{0,\pi}[e^{-\lambda \underline{g}_\infty}]=0 \quad \quad \forall \lambda >0
\end{equation*}
which implies that $\underline{g}_\infty=\infty$ $\:\bb{P}_{0,\pi}$-a.s. Hence, by Proposition \ref{Main Proposition}, we conclude that $\xi_t$ drifts to $-\infty$ $\:\bb{P}_{0,\pi}$-a.s.\\
\bf{(ii)} Assume \eqref{MainThm(ii)} holds. Then, by the proof of part $(i)$, $\xi$ drifts to $-\infty$ under $\Hat{\bb{P}}_{0,\pi}$. Hence, $\xi$ drifts to $+\infty$ $\:\bb{P}_{0,\pi}$-a.s. \\
Now, suppose that $\xi$ drifts to $+\infty$ $\bb{P}_{0,\pi}$-a.s. Then, $\xi$ drifts to $-\infty$ under $\Hat{\bb{P}}_{0,\pi}$-a.s, and by Proposition \ref{Main Proposition}, we have $\Hat{\bb{P}}_{0,\pi}(\bar{g}_\infty < \infty)=1$ which implies \eqref{MainThm(ii)} using an analogue of \eqref{integral1} obtained by using \eqref{last time at supremum of dual}.  \\
\bf{(iii)} Assume \eqref{MainThm(iii)} holds. First equality implies that $\bar{g}_\infty =\infty$ $\:\bb{P}_{0,\pi}$-a.s. by (\ref{integral1}). Similarly, for the process $((\xi,\Theta),\Hat{\bb{P}})$, second equality implies that $\bar{g}_\infty = \infty$ $\:\Hat{\bb{P}}_{0,\pi}$-a.s. which implies that $\underline{g}_\infty= \infty$ $\:\bb{P}_{0,\pi}$-a.s. Then, by part $(iii)$ of Proposition \ref{Main Proposition}, $\:\xi$ oscillates. The necessity part follows similar as part $(i)$ and $(ii)$. $\hfill \Box$

For Theorem \ref{Main Theorem}, one of our assumptions is that the reversibility property holds for the MAP. Now, we replace this property with the weak reversibility property.\\
\bf{Assumption 4*} Assume that there exists another MAP $((\xi,\Theta),\tilde{\bb{P}})$ such that Assumption 1,2,3 hold and it satisfies the weak reversibility property (\ref{WeakReversibility}) with $((\xi,\Theta),\bb{P})$.

  By the first part of Lemma \ref{Duality Lemma}, it can be deduced that\begin{equation}\label{duality}
    \bb{E}_{0,\pi}\left[e^{-\lambda \underline{g}_{e_q}}\right]=\widetilde{\bb{E}}_{0,\pi}\left[e^{-\lambda(e_q-{\Bar{g}}_{e_q})}\right].
\end{equation} 
Indeed, let $\hat{\bb{P}}_{0,\pi}$ denotes the law of $(-\xi,\Theta)$ under $\tilde{\bb{P}}_{0,\pi}$. Then, the second part of Lemma \ref{Duality Lemma} implies that $\underline{g}_t$ under $\bb{P}_{0,\theta}$ is equal in distribution to $t-{\underline{g}}_t$ under $\bb{\hat{P}}_{0,\theta}$ which is equal in distribution to $t-{\bar{g}}_t$ under $\tilde{\bb{P}}_{0,\theta}$. Then, by \eqref{Kyprianou Proposition} and (\ref{duality}), we have the following identities for the process $((\xi,\Theta),\bb{P})$\begin{align*}
    \bb{E}_{0,\pi}[e^{-\lambda \Bar{g}_{e_q}}]&=&\int_{\bb{R}^+\times S \times \bb{R}^+}V_\pi^+(dr,dv,dz)e^{-\lambda r -qr}[q\ell^+(v)+\eta_v^+(\bf{1}_{\{e_q<\zeta^+\}})] \\
    \bb{E}_{0,\pi}[e^{-\lambda \underline{g}_{e_q}}]&=&\int_{\bb{R}^+\times S \times \bb{R}^+}\widetilde{V}_\pi^+(dr,dv,dz)e^{ -qr}[q\widetilde{\ell}^+(v)+\widetilde{\eta}_v^+(e^{-\lambda e_q} \bf{1}_{\{e_q<\tilde{\zeta}^+\}})]
\end{align*}
where\begin{equation*}
    \widetilde{V}_\pi^+(dr,dv,dz)= \widetilde{\bb{E}}_{0,\pi}\left[ \int_0^{{\Bar{L}}_\infty} \bf{1}_{\{ {\Bar{L}}^{-1}_s\in dr \, , \, \Theta^+_s\in dv \, , \, \xi^+_s\in dz \} }ds \right].
\end{equation*}
Similarly, for the process $((\xi,\Theta),\Check{\bb{P}})$ where $\Check{\bb{P}}_{0,\theta}$ is the law of $(-\xi,\Theta)$ under ${{\bb{P}}}_{0,\theta}$, we have the following identities as\begin{align*}
    \Check{\bb{E}}_{0,\pi}[e^{-\lambda {\Bar{g}}_{e_q}}]&=&\int_{\bb{R}^+\times S \times \bb{R}^+}{V}_\pi^-(dr,dv,dz)e^{-\lambda r -qr}[q {\ell}^-(v)+ {\eta}_v^-(\bf{1}_{\{e_q<\zeta^-\}})]
\\
        \Check{\bb{E}}_{0,\pi}[e^{-\lambda {\underline{g}}_{e_q}}]&=&\int_{\bb{R}^+\times S \times \bb{R}^+}\tilde{V}_\pi^-(dr,dv,dz)e^{ -qr}[q \tilde{\ell}^-(v)+ \tilde{\eta}_v^-(e^{-\lambda e_q} \bf{1}_{\{e_q<\tilde{\zeta}^-\}})]
\end{align*}
where \begin{equation*}
\widetilde{V}_\pi^-(dr,dv,dz)= \widetilde{\bb{E}}_{0,\pi}\left[ \int_0^{{\underline{L}}_\infty} \bf{1}_{\{ {\underline{L}}^{-1}_s\in dr \, , \, \Theta^-_s\in dv \, , \, \xi^-_s\in dz \} }ds \right].
\end{equation*}
The following theorem is a generalized version of Theorem \ref{Main Theorem} since Assumption 4* is a weaker condition than Assumption 4.
\begin{Thm}\label{Generalised Theorem}
Let $((\xi,\Theta),\bb{P})$ be a MAP such that Assumptions 1,2,3, and 4* satisfied. Then, 
\begin{enumerate}[label=(\roman*)]
    \item $\xi_t$ drifts to $-\infty$ $\:\bb{P}_{0,\pi}$-a.s.,  if and only if\begin{equation*}
    \int_{S\times \bb{R}^+}U_\pi^+(dv,dz)\eta^+_v(\zeta^+ =+\infty)=1  \quad \text{and} \quad  \int_{S\times \bb{R}^+}\widetilde{U}_\pi^+(dv,dz)\widetilde{\eta}^+_v(\tilde{\zeta}^+ =+\infty)=1
\end{equation*}
\item  $\xi_t$ drifts to $+\infty$ $\:\bb{P}_{0,\pi}$-a.s.,  if and only if\begin{equation*}
   \int_{S\times \bb{R}^+}U_\pi^-(dv,dz)\eta^-_v(\zeta^- =+\infty)=1 \quad \text{and} \quad  \int_{S\times \bb{R}^-}\widetilde{U}_\pi^-(dv,dz)\widetilde{\eta}^-_v(\tilde{\zeta}^- =+\infty)=1
\end{equation*}
\item $\xi_t$ oscillates $\bb{P}_{0,\pi}$-a.s., if and only if\begin{equation*}
    \int_{S\times \bb{R}^+}U_\pi^+(dv,dz)\eta^+_v(\zeta^+ =+\infty)=0 \quad \text{and} \quad \int_{S\times \bb{R}^+}{U}_\pi^-(dv,dz){\eta}^-_v(\zeta^- =+\infty)=0
\end{equation*}
\end{enumerate}
\end{Thm}
\bf{Proof:} The necessity parts follow from Proposition $\ref{Main Proposition}$ in view of the following proof of the sufficiency parts.\\
\bf{(i)} As shown in \eqref{integral1} in the proof of the first part of Theorem \ref{Main Theorem}, the hypothesis \begin{equation*}
    \int_{S\times \bb{R}^+}U_\pi^+(dv,dz)\eta^+_v(\zeta^+ =+\infty)=1
\end{equation*}
implies that \begin{equation*}
    \lim_{\lambda \to 0^+} \bb{E}_{0,\pi}[e^{-\lambda \Bar{g}_\infty}]=\int_{S \times \bb{R}^+}U_\pi^+(dv,dz)\eta_v^+(\zeta^+=+\infty)=1.
\end{equation*}
It implies that $\Bar{g}_\infty$ is $\bb{P}_{0,\pi}$-a.s. finite. Similar arguments can be applied to the MAP $((\xi,\Theta),\tilde{\bb{P}})$ to get that\begin{equation*}
      \lim_{\lambda \to 0^+} \tilde{\bb{E}}_{0,\pi}[e^{-\lambda {\Bar{g}}_\infty}]=\int_{S\times \bb{R}^+}\widetilde{U}_\pi^+(dv,dz)\widetilde{\eta}^+_v(\tilde{\zeta}^+ =+\infty)=1
\end{equation*}
Hence, ${\Bar{g}}_\infty$ is finite $\tilde{\bb{P}}_{0,\pi}$-a.s. Then, as in the proof of the first part of Theorem \ref{Main Theorem}, we can conclude that
\begin{equation*}
    \bb{E}_{0,\pi}[e^{-\lambda \underline{g}_{\infty}}]=\lim_{q\to 0} \tilde{\bb{E}}_{0,\pi}[e^{-\lambda(e_q-\bar{g}_{e_q})}]=0
\end{equation*}
Hence, $\underline{g}_\infty$ is infinite $\bb{P}_{0,\pi}$-a.s.  By Proposition \ref{Main Proposition}, $\xi_t$ drifts to $-\infty$ $\bb{P}_{0,\pi}$-a.s.\\
\bf{(ii)} By the proof of the first part, we have that $\xi_t$ drifts to $-\infty$ under $\Check{\bb{P}}_{0,\pi}$; hence, $\xi_t$ drifts to $+\infty$ $\:\bb{P}_{0,\pi}$-a.s.\\
\bf{(iii)} Same as part $(iii)$ of Theorem \ref{Main Theorem}. $\hfill$ $\Box$
\section{Examples}
Let $\mathcal{H}$ be a locally compact subset of $\bb{R}^d\backslash\{ 0\}$ where $d\geq 1$. Consider an $\mathcal{H}$-valued c\'{a}dl\'{a}g Markov process $(X,P) = ((X_t)_{t\geq 0}, ( P_x)_{x\in\mathcal{H}})) $ killed at the origin with $P_x(X_0=x) = 1$. $X$ is called an \it{ $\mathcal{H}$-valued self-similar Markov process} (ssMp for short) with index $\alpha>0$ if for all $c,t>0$ and $x\in \mathcal{H}$, we have the following scaling property\begin{equation*}
    (X_t,P_x) = (cX_{c^{-\alpha}t}, P_{c^{-1}x})
\end{equation*}
The scaling property implies that $\mathcal{H}=c\mathcal{H}$ for every $c>0$; hence, $\mathcal{H}$ is necessarily a cone of $\bb{R}^d\backslash\{0\}$ in the form of \begin{equation*}
    \mathcal{H}=\Phi(\bb{R}\times \mathcal{S})
\end{equation*}
where $S$ is a locally compact subset of unit sphere $S^{d-1}$ in $\bb{R}^d$ and $\Phi$ is a homeomorphism from $\bb{R}\times S^{d-1} $ to $\bb{R}\backslash \{ 0 \}$ defined by $\Phi(y,\theta)=\theta e^y$.

One way to examine $\mathcal{H}$-valued self-similar Markov processes is to use skew-product type decomposition of the process. Consequently, the radial distance of the process becomes a real-valued stochastic process modulated by the angular part of the process. This approach forms the basis of Lamperti-Kiu transform of the self-similar Markov processes. The origin of Lamperti-Kiu transform is due to \cite{lamperti1972semi} for the case $d=1$, and a higher dimensional version can be found in \cite{kiu1980semi}. We give the generalized version of this transform in Theorem \ref{Lamperti-Kiu} below, restating \cite[Thm.2.3]{alili2017inversion} for the sake of completeness, to reveal the correspondence between $\bb{R}\times S$-valued MAPs, where $S$ is a locally compact subset of the unit sphere $S^{d-1}$, and self-similar Markov processes in terms of the homeomorphism $\Phi$ and the skew-product decomposition.
\begin{Thm}\label{Lamperti-Kiu}
\cite[Thm.2.3] {alili2017inversion} Let $\alpha >0$ and $((\xi,\Theta),\bb{P})$ be a MAP in $S\times \bb{R}$ with an extra isolated state $\partial$ and lifetime $\zeta$. Define a process $X$ by \begin{equation*}
    X_t=
 \begin{cases}
        \Theta_{\tau_t}e^{\xi_{\tau_t}} & \text{if  }  t< \int_0^\zeta e^{\alpha \xi_s}ds  \\
        0 &  \text{if   }  t\geq \int_0^\zeta e^{\alpha \xi_s}ds
    \end{cases}
    \end{equation*}
where the time change is given by\begin{equation*} \tau_t := \inf \{ s \, : \, \int_0^s e^{\alpha \xi_u}du>t\}
\end{equation*}
for $t< \int_0^\zeta e^{\alpha \xi_s}ds$. Define $P_x = \bb{P}_{x/||x||,\log ||x||}$ for $x\in \mathcal{H}$ and $P_0 = \bb{P}_{\partial}$. Then, the process $(X,P_x)$ is an $\mathcal{H}$-valued ssMp with index $\alpha>0$ and lifetime $\int_0^\zeta e^{\alpha \xi_s}ds$.\\
Conversely, let $(X,P_x)$ be an $\mathcal{H}$-valued self-similar Markov process with index $\alpha>0$ and lifetime $\zeta'$. Define the process $(\xi,\Theta)$ by \begin{equation*}
     \quad 
 \begin{cases}
        \xi_t = \log ||X_{A_t}|| \quad \text{and} \quad \Theta_t = \frac{X_{A_t}}{||X_{A_t}||} & \text{if  }  t< \int_0^{\zeta'} \frac{ds}{||X_s||^\alpha}  \\
        (\xi_t,\Theta_t)=\partial &  \text{if   }  t\geq \int_0^{\zeta'} \frac{ds}{||X_s||^\alpha}
    \end{cases}
\end{equation*}
where the time change is given by\begin{equation*}
    A_t := \inf \{  s \, : \, \int_0^s \frac{du}{||X_u||^\alpha}>t\} 
\end{equation*} for $t< \int_0^{\zeta'} e^{\alpha \xi_s}ds$. Define $\bb{P}_{y,\theta}= P_{\theta e^y}$ for $\theta \in S$ and $y\in \bb{R}$, and $\bb{P}_\partial = P_0$. Then, the process $((\xi,\Theta),\bb{P}_{y,\theta})$ is a MAP in $S\times \bb{R}$ with lifetime $ \int_0^{\zeta'} \frac{ds}{||X_s||^\alpha}$.
\end{Thm}
Note that whether the lifetime $\zeta'$ of the ssMp $X$ is finite or not, either $ \int_0^{\zeta'} \frac{ds}{||X_s||^\alpha}=\infty$ and the corresponding MAP $(\xi,\Theta)$ has infinite lifetime, or $ \int_0^{\zeta'} \frac{ds}{||X_s||^\alpha}<\infty$ and the corresponding MAP $(\xi,\Theta)$ has finite lifetime $\zeta = \int_0^{\zeta'} \frac{ds}{||X_s||^\alpha}$ \cite{alili2017inversion}.

An example of a well-known class of self-similar Markov processes $(X,P)$ and the corresponding MAP $((\xi,\Theta),\bb{P})$ via Theorem \ref{Lamperti-Kiu} is given as follows. We say that $((\xi,\Theta),\bb{P})$ and $(X,P)$ have \it{skew-product property} if the transition probabilities of $(\xi,\Theta)$ are of the form
\begin{eqnarray}\label{skew}
    \bb{P}_{y,\theta}(\xi_t\in dx, \Theta_t\in d\omega) &=& e^{-\lambda t}\, \bb{P}_y^{\xi'}(\xi'_t\in dx)\bb{P}_\theta^{\Theta'}(\Theta'_t\in d\omega) \nonumber \\
    \bb{P}_{y,\theta}((\xi_t,\Theta_t)=\partial) &=& 1-e^{-\lambda t}
\end{eqnarray} 
where $\lambda \geq 0$ is some constant, $(\xi',\bb{P}^{\xi'}_y)$ is a non-killed real Lévy process, and $(\Theta',\bb{P}^{\Theta'}_\theta)$ is a Markov process on $S$ with infinite lifetime. One of the important properties of MAPs with skew-product property is that they have the reversibility property (\ref{Reversibility}) due to \cite[Prop.3.2]{alili2017inversion}.

In this section, we will give two important examples of a ssMp with skew-product property and show that their corresponding MAPs satisfy the assumptions of our main result Theorem \ref{Main Theorem}.

\subsection{Free 2-dimensional Bessel Processes}
Let $X = ((X_t^{(1)},X_t^{(2)}),t\geq 0 )$ be a free 2-dimensional Bessel processes where $X^{(i)}$'s are independent Bessel processes of dimension $3$. Note that $X$ is an $\mathcal{H}$-valued ssMp with index $2$ where $\mathcal{H}= [0,\infty)^2$. Consider the MAP $((\xi,\Theta),\bb{P})$ that corresponds to the free Bessel process $(X,P_x)$ via Theorem \ref{Lamperti-Kiu}. Note that the state space of the modulator $\mathcal{H}_0$ is the restriction of the unit circle to the first quadrant, that is, the set $\mathcal{H}_0:=\{(x,y)\,:\, x,y\geq 0 \, , \, x^2+y^2=1\}$. As a particular case of \cite[Remark 4.6]{alili2017inversion}, there exists a 2-dimensional Brownian motion $W=(W^{(1)},W^{(2)})$ such that\begin{equation*}
 \xi_t = \int_0^t \Theta_s^{(1)}dW_s^{(1)} + \int_0^t \Theta_s^{(2)}dW_s^{(2)}+2t   
\end{equation*}
and $\Theta_t = (\Theta^{(1)}_t,\Theta^{(2)}_t)$ satisfies the following SDE system \begin{equation}\label{SDEBessel}
    d\Theta_t = b(\Theta_t)dt + \sigma(\Theta_t)dW_t
\end{equation}
with \begin{equation}\label{CoefBessel}
    b(x,y) = \begin{bmatrix} \frac{1}{x}-\frac{5}{2}x \\ \frac{1}{y}-\frac{5}{2}y \end{bmatrix} \quad \sigma(x,y) = \begin{bmatrix}
        y^2 & -xy \\
        -xy  & x^2
    \end{bmatrix}
\end{equation}
for $(x,y)\in\mathcal{H}_0$.
\begin{Lem}\label{Bessel is a T process lemma}
    The diffusion process $(\Theta,\bb{P}_\theta)$ with SDE representation (\ref{SDEBessel}) is a T-process in the sense of \cite{meyn1993astability}.
\end{Lem}
\bf{Proof:} According to \cite[Prop 2.2]{stramer1997existence}, it is sufficient to show that $\Theta$ is irreducible with respect to Lebesgue measure, and it is a Feller process. We can deduce from \cite[Ch.XI]{revuz2013continuous} that the free Bessel process $(X,P_x)$ is Feller. Therefore, the corresponding MAP $(\xi,\Theta)$, and $\Theta$ are also Feller by \cite{kyprianou2018stable}. \\
Consider $B:=\{ (x,y):x\in (a,b),y=\sqrt{1-x^2}\}$ for some arbitrary $0\leq a<b \leq 1$. We need to show that \begin{equation}\label{43}
    \int_0^\infty \bb{P}_\theta(\Theta_t\in B)dt>0
\end{equation}
for each $\theta=(\theta_1,\theta_2)\in\mathcal{H}_0$. Note that a Bessel process of dimension at least 2 never hits $0$ \cite[\S4.35]{rogers2000diffusions}. Then, by Lamperti-Kiu transform, we have\begin{equation*}
    \{ \Theta_{\tau_t}\in B \} = \{ X_t=(X_t^{(1)},X_t^{(2)})\in C \}
\end{equation*}
where $\tau_t$ is defined in Theorem \ref{Lamperti-Kiu}, and\begin{equation}\label{44}
    C:= \{ (\alpha s,\alpha \sqrt{1-s^2}) \, | \, \alpha \in \bb{R}\backslash\{ 0 \} \, , s\in (a,b) \}
\end{equation}
Since $\tau_t$ is a continuous time change, (\ref{43}) is implied by \begin{equation}\label{45}
    \int_0^\infty P_q(X_t\in C)dt >0
\end{equation}
where $q=(q_1,q_1):=(c\theta_1,c\theta_2)$ for some $c>0$.\\
We introduce a transformation to polar coordinates given by $R_t = ||X_t||$, and $\Phi_t= \arctan(X^{(2)}_t/X^{(1)}_t)$. Then, using the fact that $X^{(1)}$ and $X^{(2)}$ are independent Bessel processes of dimension $3$, and the probability transition functions of each $X^{(i)}, i=1,2$, from \cite[Ch.8]{cinlar2011probability}, we find that
\begin{eqnarray*}
    \int_0^\infty P_q (X_t \in C)dt &=& \int_0^\infty dt \int_I d\varphi\int_{(0,\infty)}   dr\, r \, \frac{r\cos\varphi}{q_1}\left[ \frac{e^{-(r\cos\varphi-q_1)^2/2t}}{\sqrt{2\pi t}}-\frac{e^{-(r\cos\varphi+q_1)^2/2t}}{\sqrt{2\pi t}}\right]\nonumber \\
    &\quad&\quad\quad\quad\quad\quad\quad\quad \cdot\frac{r\sin\varphi}{q_2}\left[ \frac{e^{-(r\sin\varphi-q_2)^2/2t}}{\sqrt{2\pi t}}-\frac{e^{-(r\sin\varphi+q_2)^2/2t}}{\sqrt{2\pi t}}\right] 
\end{eqnarray*}
for the case $q_1,q_2\neq 0$ where $I:=(\arccos a,\arccos b)$. Note that $q_1$ and $q_2$ cannot be simultaneously zero since we take $||q||>0$. In the case that one of them is $0$, say $q_1=0$, we have 
\begin{eqnarray*}
    \int_0^\infty P_q (X_t \in C)dt = \int_0^\infty dt \int_I d\varphi \int_{(0,\infty)}   dr\, r \, \frac{2r^2 \cos^2\varphi e^{-r^2\cos^2\varphi/2t}}{\sqrt{2\pi t^3}} \\
    \quad\quad\quad\quad\quad\quad\quad\quad \cdot\frac{r\sin\varphi}{q_2}\left[ \frac{e^{-(r\sin\varphi-q_2)^2/2t}}{\sqrt{2\pi t}}-\frac{e^{-(r\sin\varphi+q_2)^2/2t}}{\sqrt{2\pi t}}\right] 
\end{eqnarray*}
Since the integrand in both cases is strictly positive, (\ref{45}) holds. It follows that $\Theta$ is irreducible with respect to Lebesgue measure.
 $\hfill \Box$
\begin{Prop}\label{Bessel Proposition}
    The MAP $((\xi,\Theta),\bb{P})$ corresponding to free 2-dimensional Bessel process $X = ( (X_t^{(1)},X_t^{(2)}),t\geq 0 )$ where each $X^{(i)}$ is a Bessel process of dimension 3 satisfies the assumptions of Theorem \ref{Main Theorem}.
\end{Prop}
\bf{Proof:} According to \cite{alili2017inversion}, $(X,P_x)$ satisfies the skew-product property. Hence, we first show that $\bb{P}_{0,\theta}(\zeta = \infty ) =1$ for all $\theta\in S$ to deduce that $\lambda=0$ in the expression (\ref{skew}) so that $\xi$ and $\Theta$ are independent in the corresponding MAP $((\xi,\Theta),\bb{P}_{0,\theta})$. By Theorem \ref{Lamperti-Kiu} and the fact that $(X,P_x)$ has infinite lifetime, $(\xi,\Theta)$ has infinite lifetime if the following holds\begin{equation*}
    \int_0^\infty \frac{ds}{||X_s||^2}=\infty.
\end{equation*}
According to \cite[\S 11.1]{revuz2013continuous}, we have the following law of the iterated logarithm\begin{equation}\label{47}
    \limsup_{t\to \infty}\frac{||X_t||}{\sqrt{2t \ln\ln t}}=1.
\end{equation}
Hence, for each $\epsilon>0$, $\exists t'>0$ such that $\forall t>t'$, we have\begin{equation*}
    ||X_t||^2<2t\ln\ln t + \epsilon.
\end{equation*}
Choose $\epsilon<e$. Then,\begin{eqnarray}
    \int_0^\infty \frac{ds}{||X_s||^2} &\geq& \int_{e}^\infty \frac{ds}{||X_s||^2}\nonumber\\
    &\geq& \int_{e}^\infty \frac{ds}{2s\ln\ln s +\epsilon} \nonumber \\
    &\geq& \int_{e}^\infty \frac{ds}{2s\ln\ln s + s} \nonumber \\
    &=& \int_{1}^\infty \frac{du}{2\ln u +1}\nonumber \\
    &\geq& \int_{(1,\infty)} \frac{du}{2u-1}=\infty\nonumber
    \end{eqnarray}
where we use the change of variable $\ln s =u$, and the upper bound $\ln x \leq x-1$.

By the proof of Lemma \ref{Bessel is a T process lemma}, the MAP $(\xi,\Theta)$ is a Feller process. Hence, $\xi$ is quasi-left continuous, and $\Theta$ is a Hunt process \cite{getoor1968potential}. Moreover, $\Theta$ has infinite lifetime. We will prove that $\Theta$ is positive Harris recurrent to show the existence of an invariant probability measure $\pi$. Let $\mathcal{L}$ denote the infinitesimal generator of $\Theta$. Actually, using the SDE representation (\ref{SDEBessel}), we can write the explicit form of $\mathcal{L}$ such that \begin{equation*}
    \mathcal{L} = \sum_{i=1}^2 b_i(x,y)\partial_i  + \frac{1}{2}\sum_{i,j=1}^2 a_{ij}(x,y)\partial_i\partial_j 
\end{equation*}
where $a=\sigma\sigma^T$ and, $b$ and $\sigma$ are given in (\ref{CoefBessel}). It is straightforward to show that $a=\sigma$ in our case. According to \cite[Thm.4.2]{meyn1993bstability}, if there exist $c,d>0$, a positive function $V:\mathcal{H}_0\to \bb{R}^0$, and a petite set $C\subseteq \mathcal{H}_0$ with $V$ bounded on $C$ such that \begin{equation}\label{50}
    \mathcal{L}V(x,y) \leq -c + d\bf{1}_C(x,y)
\end{equation}
then, $\Theta$ is positive Harris recurrent. We have shown that $\Theta$ is a T-process in Lemma \ref{Bessel is a T process lemma}; hence, by \cite[Thm.4.4]{meyn1993bstability}, we can take $C$ as any compact set. Let $V(x,y)=x^3+y^3$. Then, we have\begin{equation*}
    \mathcal{L}V(x,y)=\frac{1}{2}(x+y)(-21+54xy)
\end{equation*}
Note that $\mathcal{L}V(x,y)$ and $V(x,y)$ are both bounded on $\mathcal{H}_0$, so (\ref{50}) is straightforward by taking $C=\mathcal{H}_0$. Hence, $\Theta$ has an invariant probability measure $\pi$, so Assumption 1 holds for $(\xi,\theta)$.

By the regularity properties of the Wiener processes, we can deduce that the process $X$ is sphere-exterior and sphere-interior regular. In other words, if $\tau_r^+:= \inf \{ t>0 \, | \, ||X_t||>r\}$ and $\tau_r^-:= \inf \{ t>0 \, | \, ||X_t||<r \}$, then, we have $P_x(\tau_1^+=0)=P_x(\tau_1^-=0)=1$ for each $x\in \mathcal{H}_0$. Hence, by Lamperti-Kiu transform, the MAP $(\xi,\Theta)$ is both upwards and downwards regular, so Assumption 2 is satisfied. 

Note that since $X^{(1)}$ and $X^{(2)}$ are independent Bessel processes of dimension 3, $||X||$ is a Bessel process of dimension 6, and Assumption 3 follows from the following maximal inequality for Bessel process\begin{equation*}
    \bb{E}\left[\sup_{t\in(0,\tau)}||X_t||\right]\leq \gamma\, \bb{E}[X_\tau]
\end{equation*} 
where $\gamma$ is a constant depending on the dimension of the Bessel process and $\tau$ is a stopping time \cite{dubins1994optimal}. Indeed, we use the definition $\xi_t= \log ||X_{A_t}||$ from Theorem \ref{Lamperti-Kiu} and the fact that $A_t$ is a continuous time change to deduce Assumption 3. Finally, the MAP has the reversibility property by \cite[Ex.B]{alili2017inversion}, so Assumption 4 holds.  $\hfill \Box$

\subsection{Radial Dunkl Processes}
We recall some definitions and properties about Dunkl processes from \cite{chybiryakov2008skew,alili2017inversion,demni2009radial}. For $\alpha\in \bb{R}^n \backslash \{0\}$, let $H_\alpha := \{ \beta\in \bb{R}^n | \alpha \cdot \beta =0 \}$ and $\sigma_\alpha$ denotes the reflection with respect to the hyperplane $H_\alpha$. A finite set $R\subset \bb{R}^n \backslash \{0\}$ is called a \it{root system} if $R\cap \bb{R}\alpha = \{ \pm \alpha \}$ $\forall \alpha \in R$ and $\sigma_\alpha(R)=R$ $\forall \alpha \in R$. The subgroup $W\subset O(\bb{R}^n)$, which is generated by the reflections $\sigma_\alpha$ for $\alpha\in R$, is called the Weyl reflection group associated with $R$, where $O(\bb{R}^n)$ is the orthogonal group of $\bb{R}^n$. Let $R^+$ be a positive subsystem of $R$ and $k$ be a non-negative function on $R$ such that $k\circ w = k$ $\forall w\in W $. Then, a Dunkl process $X^k$ is a c\'{a}dl\'{a}g Markov process with generator $L_k:= \frac{1}{2} \Delta_k$, where $\Delta_k = \sum_{i=1}^d T_i^2$ is the Dunkl Laplacian and $T_if(x):= \partial_if(x) + \sum_{\alpha\in R^+}k(\alpha)\alpha_i\frac{f(x)-f(\sigma_\alpha x)}{\alpha \cdot x}$. 

Let $\mathcal{C}$ be a connected component of $\bb{R}^n \backslash \cup_{\alpha\in R}H_\alpha$, and $\bb{R}^n / W$ be the space of $W$-orbits in $\bb{R}^n$. Let $\pi_1: \bb{R}^n \to \bb{R}^n / W$ and $\pi_2: \bb{R}^n /W \to \Bar{\mathcal{C}}$ be the canonical projections. Then, $X^W := \pi_2\circ \pi_1 (X^k)$ is a Feller self-similar process with index $2$ called as a \it{radial Dunkl process}. Let $W_t$ be a n-dimensional Wiener process. Then, $X^W$ is the unique strong solution of the following SDE\begin{equation}\label{SDEDunkl}
    d X_t^W = dW_t + \sum_{\alpha \in R_+} \frac{k(\alpha)}{<\alpha , X^W_t>}\alpha\, dt 
\end{equation}
\cite{demni2009radial,chybiryakov2006processus}. Furthermore, if $k(\alpha)\geq 1/2$ for any $\alpha\in R$, then, $P_y(X_t\in \bar{\mathcal{C}}\backslash \mathcal{C})=0$ for all $y\in \bar{\mathcal{C}}$ and $t>0$.

Throughout this section, we assume $k(\alpha)=k$ for some $k\geq 1/2$ for all $\alpha\in R_+$. We consider some important examples of radial Dunkl processes. For $n=1$, we can choose $R:=\{\pm 1\}$ as a root system. Let $R_+=\{ +1 \}$ and $\mathcal{C}:=(0,\infty)$. Then, the corresponding radial Dunkl process is the Bessel process with dimension $k$ \cite{revuz2013continuous}. In this case, (\ref{SDEDunkl}) takes the form\begin{equation*}
    d X^W_t = dW_t + \frac{k}{X^W_t}dt
\end{equation*}
Let $n>1$ and $(e_i)_{1\leq i \leq n}$ be the standard basis of $\bb{R}^n$. The root system $R:=\{ \pm (e_i-e_j) \, | \,  i,j\in \{1,2,\ldots,n\}, i\neq j\}$ is called $A_{n-1}$ \it{type root system}. Take $R_+ = \{ e_i-e_j\, |\, i<j \} $ and $\mathcal{C}^A:= \{ x\in \bb{R}^d \, | \, x_1 > x_2 > \ldots > x_n \}$. The resulting radial Dunkl process $X^W$ is Dyson's Brownian motion with parameter $k$ \cite{chybiryakov2008skew}.\\
The following three different types of root systems are related to Wischart processes \cite{demni2010beta,humphreys1990reflection}. \begin{itemize}
    \item The root system $R := \{ \pm e_i\, | \, i\geq 1\} \cup \{ \pm (e_i \pm e_j) \, | \, i<j \}$ is called $B_n$ \it{type root system}. Take $R_+= \{  e_i\, | \, i\geq 1\} \cup \{  (e_i \pm e_j) \, | \, i<j \}$ and $\mathcal{C}^B:=  \{ x\in \bb{R}^d \, | \, x_1 > x_2 > \ldots > x_n >0 \}$.  
    \item The root system $R := \{ \pm2 e_i\, | \, i\geq 1\} \cup \{ \pm (e_i \pm e_j) \, | \, i<j \}$ is called $C_n$ \it{type root system}. Take $R_+= \{  2e_i\, | \, i\geq 1\} \cup \{  (e_i \pm e_j) \, | \, i<j \}$ and $\mathcal{C}^C:=  \{ x\in \bb{R}^d \, | \, x_1 > x_2 > \ldots > x_n >0 \}$. 
    \item The root system $R :=  \{ \pm (e_i \pm e_j) \, | \, i<j \}$ is called $D_n$ \it{type root system}. Take $R_+= \{  (e_i \pm e_j) \, | \, i<j \}$ and $\mathcal{C}^D:=  \{ x\in \bb{R}^d \, | \, x_1 > x_2 > \ldots > |x_n| >0 \}$. 
\end{itemize}
\indent
The following proposition gives the characterization of the MAP $((\xi,\Theta),\bb{P})$ corresponding to a radial Dunkl process using the SDE representation (\ref{SDEDunkl}).
\begin{Prop}\label{MAP representation of Dunkl Proposition}
Let $(\xi,\Theta)$ be the MAP corresponding to the radial Dunkl process $X^W$ via the Lamperti-Kiu transform. Then, there exists an $n$ dimensional Wiener process $W_t=(W_t^{(1)},\ldots, W_t^{(n)})$ such that $\xi$ satisfies the SDE \begin{equation*}
    d\xi_t= \sum_{i=1}^n \Theta_t^{(i)}dW_t^{(i)} + \left(\sum_{i=1}^n (\Theta_t^{(i)} \sum_{\alpha\in R_+} \frac{k(\alpha)}{<\alpha,\Theta_t>}\alpha_i) + n-2\right)dt 
\end{equation*}
and $\Theta=(\Theta^{(1)},\ldots,\Theta^{(n)})$ satisfies the SDE system \begin{multline*}
    d\Theta_t^{(i)} = dW_t^{(i)} - \Theta_t^{(i)} \sum_{j=1}^n dW_t^{(j)}  \\ 
    +\left( \sum_{\alpha\in R_+} \frac{k(\alpha)}{<\alpha,\Theta_t>}\alpha_i
    - \Theta_t^{(i)}\sum_{j=1}^n (\Theta_t^{(j)} \sum_{\alpha\in R_+} \frac{k(\alpha)}{<\alpha,\Theta_t>}\alpha_j) -\frac{\Theta_t^{(i)}}{2} \right) dt
\end{multline*}
\end{Prop}
\bf{Proof:} Follows from Itô formula in view of (\ref{SDEDunkl}).$\hfill\Box$

We consider the case $d=2$ for simplicity. Let $(\xi,\Theta)$ be the corresponding MAP of the radial Dunkl process $X^W=(X^{(1)},X^{(2)})^W$ in $\bb{R}^2 \backslash \{ 0 \}$. We will show that $(\xi,\Theta)$ satisfies the assumptions of Theorem \ref{Main Theorem} for the cases in which the root system $R$ is of type $A_1,\, B_2,\, C_2$, or $D_2$.
\begin{Lem}\label{Dunkl Lemma}
    The diffusion process $(\Theta,\bb{P}_\theta)$ is a T-process in the sense of \cite{meyn1993astability}.
\end{Lem}
\bf{Proof:} The process $X^W$ is Feller by \cite{chybiryakov2008skew}. Hence, according to the proof of Lemma \ref{Bessel is a T process lemma}, it is sufficient to show that \begin{equation}\label{54}
    \int_0^\infty P_q(X_t\in C)dt>0
\end{equation}
where $q$ and $C$ are defined as in (\ref{44}) and (\ref{45}) for appropriate $a,b$ in each case of the root system $R$. Using the explicit form of the transition density function of $X_t^W$ from \cite[Eq.6]{chybiryakov2008skew}, the change of variables $R_t=||X^W_t||, \, \Phi_t=\arctan (X^{(2)},X^{(1)})$, and the fact that $P_y(X_t\in \bar{\mathcal{C}}\backslash \mathcal{C})=0$, we get\begin{multline*}
    \int_0^\infty P_q(X_t\in C)dt = \int_0^\infty dt \int_I d\varphi \int_{(0,\infty)} dr \, \frac{r}{c_k t^{\gamma +1}}\exp\left(-\frac{c^2+r^2}{2t}\right)\\
    D_k^W\left(\frac{q}{\sqrt{t}},\frac{(r\cos\varphi,r\sin\varphi)}{\sqrt{t}}\right)\omega_k((r\cos\varphi,r\sin\varphi))
\end{multline*} 
where $I=(\arccos a,\arccos b)$, $D_k(u,v)>0$ is the Dunkl kernel, $\gamma= \sum_{\alpha\in R_+}k>0$, $\omega_k(y)=\prod_{\alpha\in R_+}|\alpha \cdot y | ^{2k}$ and $c_k=\int_{\bb{R}^2}e^{-|x|^2/2}\omega_k(x)dx$. The fact that $P_y(X_t\in \bar{\mathcal{C}}\backslash \mathcal{C})=0$ implies that $P_q(X_t\in \cup_{\alpha\in R_+} H_\alpha)=0$, so $\omega_k(y)>0$ and $c_k>0$ a.s. Hence, the integrand is strictly positive, and (\ref{54}) holds, which in turn implies that $\Theta$ is irreducible. $\hfill\Box$
\begin{Prop}\label{Dunkl Proposition}
    The MAPs $((\xi,\Theta),\bb{P})$ corresponding to radial Dunkl processes $(X,P_x)$ in $\bb{R}^2\backslash \{ 0 \}$ of types $A_1$, $B_2$, $C_2$, or $D_2$ satisfy the assumptions of Theorem \ref{Main Theorem}.
\end{Prop}
\bf{Proof}: According to \cite{alili2017inversion}, $(X,P_x)$ satisfies the skew-product property. Hence, we first show that $\bb{P}_{0,\theta}(\zeta = \infty)=1$ for all $\theta \in S$ to deduce that $\xi$ and $\Theta$ are independent as in the proof of Proposition \ref{Bessel Proposition}. Note that the property
$P_y(X_t \in \Bar{\mathcal{C}}\backslash\mathcal{C})=0$ for all $y\in \Bar{\mathcal{C}}$ implies that the lifetime of $(X,P_x)$ is infinite a.s. According to this fact and Theorem \ref{Lamperti-Kiu}, the MAP $(\xi,\Theta)$ has an infinite lifetime if\begin{equation}\label{56}
    \int_0^\infty \frac{ds}{||X_s||^2}=\infty.
\end{equation}
According to \cite[Ex.C]{alili2017inversion}, $||X_s||$ is a Bessel process of dimension $2+2\gamma$ where $\gamma>0$ is defined in the proof of Lemma \ref{Dunkl Lemma}. Hence, by the law of the iterated logarithm (\ref{47}), Equation (\ref{56}) holds.

We have shown that the MAP $(\xi,\Theta)$ is a Feller process in Lemma \ref{Dunkl Lemma}. Hence, the ordinate $\xi$ is quasi-left continuous, and the modulator $\Theta$ is a Hunt process. Moreover, $\Theta$ has infinite lifetime. We will prove that $\Theta$ is positive Harris recurrent to show the existence of an invariant probability measure $\pi$. Let $\mathcal{L}$ denote the infinitesimal generator of $\Theta$. Using the particular case $n=2$ of the SDE representation in Proposition \ref{MAP representation of Dunkl Proposition}, we have 
    $\mathcal{L} = \sum_{i=1}^2b_i(x,y)\partial_i + \frac{1}{2}\sum_{i=1}^2a_{ij}\partial_i\partial_j$
where \begin{equation*}
     b(x,y) = \begin{bmatrix} y^2\mathcal{A}_1-xy\mathcal{A}_2-\frac{x}{2} \\ x^2\mathcal{A}_2 - xy\mathcal{A}_1 - \frac{y}{2} \end{bmatrix} \quad \sigma(x,y) = \begin{bmatrix}
        y^2 & -xy \\
        -xy  & x^2
    \end{bmatrix}
\end{equation*}
for $\mathcal{A}_i:= \sum_{\alpha\in R^+}\frac{k}{\alpha_1x+\alpha_2y}\alpha_i$ with $i=1,2$, and $a := \sigma\sigma^T=\sigma$. Denote the corresponding state space of $\Theta$ by $\mathcal{C}_0$ for given $\mathcal{C}$. According to \cite[Thm.4.2]{meyn1993bstability}, if there exist $c,d>0$, a positive function $V:\mathcal{C}_0\to \bb{R}^+$, and a petite set $C\subset \mathcal{C}_0$ with $V$ bounded on $C$ such that \begin{equation}\label{59}
    \mathcal{L} V(x,y) \leq -c + d \bf{1}_C(x,y)
\end{equation} 
then, $\Theta$ is positive Harris recurrent. By Lemma \ref{Dunkl Lemma}, $\Theta$ is a T-process; hence, by \cite[Thm.4.4]{meyn1993bstability}, we can take $C$ as any compact set. In each case, we take $C$ as the state space of the modulator and find appropriate $V(x,y)$ for each mentioned type of finite root system such that both $\mathcal{L}V(x,y)$ and $V(x,y)$ are bounded on $C$. Then, (\ref{59}) will follow immediately.\\
\it{Type $A_1$}: In this case, $R_+=\{ (1,-1) \}$. We have $\mathcal{A}_1= \frac{k}{x-y}$ and $\mathcal{A}_2=\frac{-k}{x-y}$. Let $V(x,y)=x^3+y^3$. Then, we obtain\begin{equation*}
    \mathcal{L}V(x,y)= \frac{3}{2}(x+y)[(4k+2)xy-1]+(x+y)(12xy-3)
\end{equation*}
$\mathcal{L}V(x,y)$ and $V(x,y)$ are both finite on $\mathcal{C}_0^A$.\\
\it{Type $B_2$}: In this case, $R_+=\{ (1,-1),(1,1),(1,0),(0,1)\}$. Then, $\mathcal{A}_1= \frac{k}{x-y}+\frac{k}{x+y}+\frac{k}{x}$ and $\mathcal{A}_2= \frac{-k}{x-y}+\frac{k}{x+y}+\frac{k}{y}$. Take $V(x,y)=x^3+y^3$. Then, we get \begin{equation*}
    \mathcal{L}V(x,y)= \frac{6k}{x+y}(8x^2y^2-1)+ \frac{-3}{2}(x+y)(x-y)^2 + (x+y)(12xy-3)
\end{equation*}
In this case, $\mathcal{C}_0^B = \{ (\cos \alpha , \sin \alpha ) : \alpha \in (\pi/4,\pi/2) \}$; hence, the term $x+y$ is always greater than 1. Therefore, $\mathcal{L}V(x,y)$ and $V(x,y)$ are both finite on $\mathcal{C}_0^B$.\\
\it{Type $C_2$}: In this case, $R_+= \{ (1,-1),(1,1),(2,0),(0,2)\}$. Then, we have $\mathcal{A}_1= \frac{k}{x-y}+\frac{k}{x+y}+\frac{2k}{x}$ and $\mathcal{A}_2= \frac{-k}{x-y}+\frac{k}{x+y}+\frac{2k}{y}$. Take $V(x,y)=x^3+y^3$. Then, we have \begin{equation*}
    \mathcal{L}V(x,y)= \frac{12k}{x+y}(6x^2y^2-1)-\frac{3}{2}(x-y)^2(x+y)+(x+y)(12xy-3)
\end{equation*}
In this case, we have $\mathcal{C}_0^C=\mathcal{C}_0^B$; hence, the term $x+y$ is greater than 1. Therefore, $\mathcal{L}V(x,y)$ and $V(x,y)$ are both finite on $\mathcal{C}_0^C$.\\
\it{Type $D_2$}: In this case, $R_+=\{ (1,-1),(1,1) \}$. We have $\mathcal{A}_1= \frac{k}{x-y}+\frac{k}{x+y}$ and $\mathcal{A}_2=\frac{-k}{x-y}$. Take $V(x,y)=(x^3+y^3)^2$. In this case, we have\begin{equation*}
    \mathcal{L}V(x,y)=48kx^2y^2(1-xy) +3(4x^2y^2-1)(1-xy)+6(x+y)^2(1-xy)(4xy-1)+ (x+y)(12xy-3)
\end{equation*}
It can be verified that $\mathcal{L}V(x,y)$ and $V(x,y)$ are both finite on $\mathcal{C}_0^D$.\\
Hence, Assumption 1 holds for MAPs $(\xi,\Theta)$ corresponding to respective root systems. 

$||X_s||$ is a Bessel process of dimension $2+\gamma \geq 3$, and $\xi$ and $\theta$ are independent. Then, Assumption 2 follows from the regularity properties of Bessel processes from \cite[Ch.5.48]{rogers2000diffusions}.
Assumption 3 follows exactly as in the proof of Proposition \ref{Bessel Proposition}. Finally, the MAP has reversibility property by \cite[Ex.C]{alili2017inversion}, so Assumption 4 holds.  $\hfill\Box$
\\
\bf{Remark 3} We have established that  $\Theta$ is positive Harris recurrent for both examples in our pursuit of verifying the assumptions of Theorem \ref{Main Theorem}. Therefore, the angular part $X_t/||X_t||$ of each self-similar Markov process is also positive Harris recurrent since the time change $\tau_t$ in Theorem \ref{Lamperti-Kiu} is continuous and Assumption 3 holds.

\bibliographystyle{abbrv}
\bibliography{references}

\end{document}